\documentclass[a4paper]{article}
\usepackage{amsmath,amssymb,fancybox,longtable,graphics,amscd,multirow,theorem}
\usepackage[dvipdfm]{graphicx}
\usepackage{caption,setspace}
\newtheorem{thm}{\textsc{Theorem}}[section]

\newtheorem{lem}{\textsc{Lemma}}[section]
\newtheorem{cor}{\textsc{Corollary}}[section]

\def\QED{$\Box$}
\def\END{$\blacksquare$}
\def\mbi#1{\boldsymbol{#1}} 

\def\Pic{\mathop{\mathrm{Pic}}\nolimits}
\def\discr{\mathop{\mathrm{discr}}\nolimits}
\def\rk{\mathop{\mathrm{rank}}\nolimits}
\def\sign{\mathop{\mathrm{sgn}}\nolimits}

\def\C#1{\mathbb{Z}\slash #1\mathbb{Z}}
\theorembodyfont{\rmfamily}
\newtheorem{remark}{\textsc{Remark}}

\theorembodyfont{\rmfamily}
\newtheorem{defn}{\textsc{Definition}}[section]
\newtheorem{eg}{Example \vspace{2mm}}[section]

\begin{document}
\title{Lattice duality for families of $K3$ surfaces associated to transpose duality}
\author{Makiko Mase}
\date{\small \begin{tabular}{l} Key Words: $K3$ surfaces, toric varieties, Picard lattices \end{tabular}}
\maketitle
\begin{abstract}
The aim of this article is to show that the transpose-dual pairs in the sense of Ebeling-Ploog of singularities $(Z_{1,0},\, Z_{1,0})$, $(U_{1,0},\, U_{1,0})$, $(Q_{17},\, Z_{2,0})$, $(W_{1,0},\, W_{1,0})$ that are concluded to be polytope-dual by the author are actually lattice-dual. 
\end{abstract}
\section{Introduction} \label{Introduction}
The transpose duality studied by Ebeling and Ploog \cite{EbelingPloog} is a duality for bimodal and other singularities analogous to Arnold's strange duality for unimodal singularities \cite{Arnold75}. 
A polytope duality appears in a study by Batyrev \cite{BatyrevMirror} to construct a mirror partner : this is a duality between Calabi-Yau hypersurfaces in toric varieties. 
Dolgachev \cite{DolgachevMirror} studies a duality between lattices that have a deep relation with $K3$ surfaces. 
Namely, primitive sublattices $M$ and $M'$ of the $K3$ lattice $\Lambda_{K3}$ are dual in this sense if an isometry 
\[
M_{\Lambda_{K3}}^\perp \simeq U\oplus M'
\]
holds, which implies $\rk M+\rk M'=20$. 
Here $U$ is the even unimodular hyperbolic lattice of rank $2$. 

In \cite{Mase16I}, it is shown that for the transpose-dual pairs 
\[
(Z_{1,0},\, Z_{1,0}),\, 
(U_{1,0},\, U_{1,0}),\, 
(Q_{17},\, Z_{2,0}),\,  
(W_{1,0},\, W_{1,0}), 
\]
one can take deformations $F,\, F'$ of their defining polynomials and reflexive polytopes $\Delta,\, \Delta'$ such that the conditions
\[
(*) \quad
\Delta^* \simeq \Delta',\quad
\Delta_F\subset \Delta,\quad 
\Delta_{F'}\subset \Delta', \quad 
\textrm{and} \quad \rk L_0(\Delta)=0
\]
hold. 
Here $\Delta^*$ denotes the polar dual polytope of $\Delta$, $\Delta_{\!F}$ is the Newton polytope of $F$, and $L_0(\Delta)$ is the cokernel of the restriction map of the $(1,1)$ Hodge component of the ambient space to a $K3$ surface. 
We denote by $\rho_{\Delta},\, \rho_{\Delta'}$ the rank of the Picard lattices $\Pic_{\Delta},\, \Pic_{\Delta'}$ of the families $\mathcal{F}_{\Delta},\, \mathcal{F}_{\Delta'}$ of $K3$ surfaces corresponding to the polytopes $\Delta,\, \Delta'$, respectively. 
Moreover, $\rho_{\Delta}+\rho_{\Delta'}=20$. 
We are interested in a duality of lattices associated to these polytope-dual pairs. 

The aim of this article is to show the following: \\

\noindent
{\bf Main Theorem} (Theorem \ref{MainThmPf})\, 
{\it Let $(B,\, B')$ be a pair of singularities}
\[
(Z_{1,0},\, Z_{1,0}),\, 
(U_{1,0},\, U_{1,0}),\, 
(Q_{17},\, Z_{2,0}),\,  \textit{or \,}
(W_{1,0},\, W_{1,0}). 
\]
{\it The families $\mathcal{F}_{\!\Delta}$ and $\mathcal{F}_{\!\Delta'}$ of $K3$ surfaces obtained in \cite{Mase16I} associated to $(B,\, B')$ are lattice dual. }

Moreover, there exist negative-definite even lattices $K_1$ and $K_2$ such that $\rk{K_1}=15$ and $\discr{K_1}=-18$, and $\rk{K_2}=16$ and $\discr{K_2}=4$ as in  Table \ref{LatticeExcPairsIntro}. 
Here, $\discr{K}$ indicates the discriminant of the lattice $K$. 
\[
\begin{array}{c @{\hspace{.5cm}} c @{\hspace{.5cm}} c @{\hspace{.5cm}} c @{\hspace{.5cm}} c @{\hspace{.5cm}} c @{\hspace{.5cm}} c} 
B & \Pic_\Delta & \rho_\Delta & |\discr| & \rho_{\Delta'} & \Pic_{\Delta'} & B' \\ 
\hline \hline
Z_{1,0} & U\oplus D_5\oplus E_7 & 14 & 8 & 6 & U\oplus A_1\oplus A_3 & Z_{1,0}\\
\hline
U_{1,0} & U\oplus K_1& 17  & 18 & 3 & \left( \begin{smallmatrix}0 & 3 \\ 3 & -2\end{smallmatrix}\right)\oplus A_1 & U_{1,0}\\
\hline
Q_{17} & U\oplus E_6\oplus E_7 & 15 & 6 & 5 & U\oplus A_1\oplus A_2 & Z_{2,0}  \\
\hline
W_{1,0} & U\oplus K_2 & 18 & 4 & 2 & \left( \begin{smallmatrix}0 & 2 \\ 2 & -2 \end{smallmatrix}\right) & W_{1,0}\\
\hline
\end{array}
\]
\begingroup
\captionof{table}{{\it Picard lattices for lattice dual pairs}}\label{LatticeExcPairsIntro}
\endgroup
\noindent

\bigskip
Section \ref{Preliminary} is devoted to summarise the result of \cite{Mase16I} and basic facts on toric geometry and lattice theory.  
The proof of the main theorem is given in Theorem \ref{MainThmPf} in section \ref{MainThm}. 
The method of the proof is by general theory of toric geometry (in particular, singularities) for computation of lattices, and by standard lattice theory for verification of the lattice being primitively embedded into the $K3$ lattice. \\

\noindent
\begin{ackn}\\
{\rm {
The author would thank to Professor N. Aoki who was reading through the first draft of the paper carefully and giving her many helpful comments, and to the anonymous referees for cautious reading together with mindful comments. 
}}
\end{ackn}
\section{Preliminary} \label{Preliminary}
A {\it lattice} is a finitely generated abelian group with a non-degenerate pairing.  

A {\it Gorenstein $K3$ surface} is a compact complex connected $2$-dimensional algebraic variety $S$ with at most $ADE$ singularities with trivial canonical bundle and $H^1(S,\mathcal{O}_S)=0$. 
A nonsingular Gorenstein $K3$ surface is called a {\it $K3$ surface}. 

Denote by $\mathbb{P}_\Delta$ the toric $3$-fold associated to a $3$-dimensional polytope $\Delta$. 
Suppose that $\Delta$ is integral, and contains the origin in its interior. 
Recall that by definition, the polytope $\Delta$ is {\it reflexive} if its polar dual $\Delta^*$ is also integral. 
It is shown in Theorem 2.2.24 \cite{BatyrevMirror} that the polytope $\Delta$ is reflexive if and only if the corresponding toric $3$-fold $\mathbb{P}_{\!\Delta}$ is Fano. 
If this is the case, general anticanonical sections of $\mathbb{P}_{\!\Delta}$ are Gorenstein $K3$ surface. 

Here we fix and list notation. 
See \cite{Kobayashi} and \cite{Mase16I} for more details. 
Let $\Delta$ be a $3$-dimensional reflexive polytope. 
In this article, we define that a fibre $Z$ of a family $\mathcal{F}_{\!\Delta}$ of $K3$ surfaces determined by the polytope $\Delta$ is {\it generic} if the following two conditions are satisfied: 
\begin{itemize}
\item[(1)] $Z$ is $\Delta$-regular. (See \cite{BatyrevMirror} for detail)
\item[(2)] The Picard group of $\widetilde{Z}$ is generated by irreducible components of the restrictions of the generators of the Picard group of $\widetilde{\mathbb{P}}_{\!\Delta}$. 

\end{itemize}
Here, $\widetilde{\mathbb{P}}_{\!\Delta}$ and $\tilde{Z}$ are respectively the minimal models of a simultaneous resolution of $\mathbb{P}_{\!\Delta}$ and of a generic member $Z$ in $\mathcal{F}_{\!\Delta}$. 
Let $r\colon H^{1,1}(\widetilde{\mathbb{P}}_{\!\Delta}, \mathbb{Z}) \to H^{1,1}(\tilde{Z}, \mathbb{Z})$ be the restriction map of the Hodge components. 

It is proved in \cite{BatyrevMirror} that $\Delta$-regularity is a general condition. 
The general sections satisfy the condition (2) (see \cite{Bruzzo-Grassi}). 
Note that all Picard lattices of the minimal models of any generic members are isometric.
\begin{itemize}
\item $\mathcal{F}_\Delta$ is the family of (Gorenstein) $K3$ surfaces 
\footnote{We use a convention that we identify a Gorenstein $K3$ surface and its minimal resolution that is a $K3$ surface. } 
parametrised by the complete anticanonical linear system $\mathcal{F}_\Delta \to |{-}K_{\mathbb{P}_\Delta}|$. 
\item $L_0(\Delta)$ is the cokernel of $r$. 
The rank of $L_0(\Delta)$ is calculated by a formula given in \cite{Kobayashi}. 
\end{itemize}

\begin{defn}
Let $\Delta$ be a $3$-dimensional reflexive polytope. \\
$(1)$\, The {\it Picard lattice $\Pic_{\Delta}$ of the family $\mathcal{F}_{\Delta}$} is the Picard lattice of the minimal model of any generic member in $\mathcal{F}_{\Delta}$. \\
$(2)$\, Denote by $\rho_{\Delta}$ the rank of $\Pic_{\Delta}$ and call it the {\it Picard number} of  $\mathcal{F}_{\Delta}$. \END
\end{defn}

Note that the Picard lattice of the family is defined independently of the choice of generic members. 
Suppose the Picard group of $\widetilde{\mathbb{P}}_{\!\Delta}$ is generated by toric divisors $\widetilde{D}_1,\ldots ,\widetilde{D}_s$. 
Then the Picard lattice $\Pic_\Delta$ of the family $\mathcal{F}_\Delta$ is generated by the  irreducible components of restricted divisors $D_1,\ldots ,D_s$ of $\widetilde{D}_1,\ldots ,\widetilde{D}_s$ by the requirement $(2)$ for generic fibres. 

In order to avoid the repetition, we define a notion of ``lattice duality''. 

\begin{defn}
For $3$-dimensional reflexive polytopes $\Delta$ and $\Delta'$, the families $\mathcal{F}_\Delta$ and $\mathcal{F}_{\Delta'}$ of $K3$ surfaces are called {\it lattice  dual} if an isometry 
\[
{(\Pic_\Delta)}_{\Lambda_{K3}}^\perp\simeq U\oplus\Pic_{\Delta'}
\]
holds. \END
\end{defn}

We say that a transpose-dual pair $(B, B')$ of singularities are {\it polytope-dual} if the singularities are compactified as general anticanonical sections $F,\, F'$ of $\mathbb{P}_{\!\Delta}$ and $\mathbb{P}_{\!\Delta^*}$, resp., where $\Delta$ is a reflexive polytope, satisfying the conditions 
\[
(*) \quad
\Delta^* \simeq \Delta',\quad
\Delta_F\subset \Delta,\quad 
\Delta_{F'}\subset \Delta', \quad 
\textrm{and} \quad \rk L_0(\Delta)=0. 
\]
Here $\Delta_f$ denotes the Newton polytope of a polynomial $f$. 

It is proved in \cite{Mase15} that some polytope-dual pairs extend to a lattice-dual pair between families of $K3$ surfaces. 
More precisely, 
\begin{thm}{\rm (\cite{Mase15})}
For each of the transpose-dual pairs of singularities $(B, B')$ listed in Table \ref{NicePairsIntro}, there exists a reflexive polytope $\Delta=\Delta_{[MU]}$ such that the polytope duality extends to a lattice duality between the families $\mathcal{F}_\Delta$ and $\mathcal{F}_{\Delta^*}$, where the Picard lattices are also given in Table \ref{NicePairsIntro}. 
Here we use the notation $C^6_8:=\left( \begin{smallmatrix} -4 & 1 \\ 1 & -2 \end{smallmatrix}\right)$, and $\Delta_{[MU]}$ is the polytope obtained in \cite{Mase16I}. 
\[
\begin{array}{c @{\hspace{0.5cm}} c @{\hspace{0.5cm}} c @{\hspace{0.5cm}} c @{\hspace{0.5cm}} c @{\hspace{0.5cm}} c @{\hspace{0.5cm}} c} 
B & \Pic_{\Delta_{[MU]}} & \rho_{\Delta_{[MU]}} & |\discr| &  \rho_{\Delta_{[MU]}^*} & \Pic_{\Delta_{[MU]}^*} & B' \\ 
\hline \hline
Q_{12} & U\oplus E_6\oplus E_8 & 16 & 3 & 4 & U\oplus A_2 & E_{18}\\
\hline
Z_{1,0} & U\oplus E_7\oplus E_8 & 17 & 2 & 3 & U\oplus A_1 & E_{19}\\
\hline
E_{20} & U\oplus E_8^{\oplus 2} & 18 & 1 & 2 & U & E_{20}  \\
\hline
Q_{2,0} & U\oplus A_6\oplus E_8 & 16 & 7 & 4 & U\oplus C^6_8 & Z_{17}\\
\hline
E_{25} & U\oplus E_7\oplus E_8 & 17 & 2 & 3 & U\oplus A_1 & Z_{19}\\
\hline
Q_{18} & U\oplus E_6\oplus E_8 & 16 & 3 & 4 & U\oplus A_2 & E_{30}\\
\hline
\end{array}
\]
\begingroup
\captionof{table}{{\it Picard lattices for lattice dual pairs} }\label{NicePairsIntro}
\endgroup 
\end{thm}

In \cite{Mase16I}, we can pick up more candidates of lattice-dual pairs: 
\begin{thm}{\rm (\cite{Mase16I})}\label{Mase16IMainThm}
Each of the following transpose-dual pairs are polytope dual: 
\[
(Z_{1,0},\, Z_{1,0}),\, 
(U_{1,0},\, U_{1,0}),\, 
(Q_{17},\, Z_{2,0}),\,  
(W_{1,0},\, W_{1,0}). 
\]
Namely, for each of these pairs, there exist deformations $F,\, F'$ of the polynomials defining the singularities, and reflexive polytopes $\Delta,\, \Delta'$ such that the conditions
\[
(*) \quad
\Delta^* \simeq \Delta',\, 
\Delta_F\subset \Delta,\,  \Delta_{F'}\subset \Delta', \quad \textit{and} \quad \rk L_0(\Delta)=0
\]
hold. 
Moreover, $\rho_{\!\Delta}+\rho_{\!\Delta'}=20$. 
\end{thm}
\begin{eg}
: {\bf \mathversion{bold}$Z_{2,0}$ and $Q_{17}$ case} (Theorem 3.1 \cite{Mase16I})
The defining polynomials of singularities $B=Z_{2,0}$ and $B'=Q_{17}$ are $f=x^5z +xy^3 + z^2,\, f'=x^5y +y^3+xz^2$, respectively. 

Take a deformation of $f$ (resp. $f'$) as $F= W^7Y+X^5Z +XY^3 + Z^2$ (resp. $F'=W^7+X^5Y +WY^3+XZ^2$) in the weighted projective space $\mathbb{P}(1,1,3,5)$ (resp. $\mathbb{P}(1,1,2,3)$). 

Let $\Delta$ and $\Delta'$ be polytopes that are respectively convex hulls of vertices $\left\{ (-1,-1,2), (0,-1,0), (1,-1,0), (1,-1,1), (1,2,-3), (0,0,-1) \right\}$, and \\
$\left\{ (-1,2,-1), (-1,-1,1), (-1,-1,-1), (6,-1,-1), (2,1,-1), (0,-1,1)\right\}$. 

We verify these polytopes satisfy the condition $(*)$ in Theorem 2.2. 

The polytope $\Delta$ contains the Newton polytope of $F$: indeed, by taking a basis $\mbi{e}_1=(-3,3,0,0)$,\, $\mbi{e}_2=(-8,0,1,1)$,\, $\mbi{e}_3=(-6,1,0,1)$ for $\mathbb{R}^3$, one can see that monomials $W^7Y$,\,$X^5Z$,\, $XY^3$,\, $Z^2$ are respectively corresponding to vertices $(0,0,-1),\, (1,-1,1),\, (1,2,-3),\, (-1,-1,2)$. 

The polytope $\Delta'$ contains the Newton polytope of $F'$: indeed, by taking a standard basis $\mbi{e}_1'=(-1,1,0,0)$,\, $\mbi{e}_2'=(-2,0,1,0)$,\, $\mbi{e}_3'=(-3,0,0,1)$ for $\mathbb{R}^3$, one can see that monomials $W^7,\, X^5Y,\, WY^3,\, XZ^2$ are respectively corresponding to vertices $(-1,-1,-1),\, (4,0,-1),\, (-1,2,-1),\, (0,-1,1)$. 

The dual polytope $(\Delta')^*$ of $\Delta'$ is a convex hull of vertices 
\[
(-1,-3,-4),\, (0,-2,-3),\, (0,1,0),\, (1,0,0),\, (0,0,1),\,  (-1,-2,-3)
\]
that is mapped to isomorphically from $\Delta$ by a transformation of $\mathbb{R}^3$ by the matrix $T := \left(\begin{smallmatrix} 1 & 1 & 1 \\ 1 & 3 & 4 \\ 1 & 2 & 3 \end{smallmatrix}\right)$, that is, $T(x, y, z)=(x', y', z')$ for $(x,y,z)\in\Delta$ and $(x',y',z')\in\Delta'$. 

Therefore, $\Delta$ and $\Delta'$ are reflexive and $\Delta^* = \Delta'$ holds. 

By the formula in \cite{Kobayashi}, one gets $\rk L_0(\Delta)=\rk L_0(\Delta^*)=0$ because for all edges in $\Delta$ satisfy $l^*(\Gamma)l^*(\Gamma^*)=0$. 

One also can compute that $\rho_{\Delta} = 18-3 = 15$, and $\rho_{\Delta'} = 8-3 = 5$ so that $\rho_{\Delta} +\rho_{\Delta'} = 20$.  

Therefore, Theorem \ref{Mase16IMainThm} holds for the pair of singularities $(Z_{2,0},\, Q_{17})$ with polynomials $F, F'$ and polytopes $\Delta, \Delta'$ presented above. \END
\end{eg}

We would like to know whether or not they are really lattice-dual pairs. 
We state the necessary facts for computation of Picard lattices by using toric geometry. 

Let $v_1,\ldots ,v_{d+3}$ be primitive one-simplices of a $3$-dimensional reflexive polytope $\Delta$. 
Each $v_i$ is associated to a {\it toric divisor} $\widetilde{D}_i$ of $\widetilde{\mathbb{P}}_{\!\Delta}$ by $\widetilde{D}_i:=\overline{{\rm Orb}(\mathbb{R}_{\geq 0}v_i)}$, where {\rm Orb} denotes the torus-action orbit. 
Note that $\mathbb{T}$ is the $3$-torus acting on $\mathbb{P}_\Delta$. 
There is an exact sequence \cite{Oda78}
\[
\begin{array}{ccccccccc}
0 & \to & M & \to & {\rm Div}_{\mathbb{T}}(\widetilde{\mathbb{P}}_{\!\Delta}) & \to & \Pic_{\widetilde{\mathbb{P}}_{\!\Delta}} & \to & 0 \\
 & & \Arrowvert & & \downarrow & & \downarrow \\
 & & M & \to & \displaystyle\bigoplus_{i=1}^{d+3}\mathbb{Z}\widetilde{D}_i & \to & A_2(\widetilde{\mathbb{P}}_{\!\Delta}) & \to & 0 
\end{array}
\]
where $M\simeq\mathbb{Z}^3$ with the standard basis $\{e_1,\, e_2,\, e_3\}$, $\Delta\subset M\otimes_{\mathbb{Z}}\mathbb{R}$, and ${\rm Div}_{\mathbb{T}}(\widetilde{\mathbb{P}}_{\!\Delta})$ is the set of all toric divisors in $\widetilde{\mathbb{P}}_{\!\Delta}$. 
The toric divisors are dependent by a linear system with three equations: 
\begin{equation}
\sum_{i=1}^{d+3}(v_i,\, e_j) \widetilde{D}_i= 0 \qquad j=1,2,3, \label{ToricLinearRelation}
\end{equation}
where $( \, , \, )$ is the standard inner product in $\mathbb{R}^3$. 
Thus, the Picard group of $\widetilde{\mathbb{P}}_{\!\Delta}$ is generated by toric divisors that are linearly independent. 
For the reflexive polytope $\Delta$ obtained in Theorem \ref{Mase16IMainThm}, it is concluded that $\rk L_0(\Delta)=0$, which means that the restriction map $r$ is surjective,  thus the Picard lattice $\Pic_\Delta$ of the family $\mathcal{F}_\Delta$ is generated by the restriction of linearly independent toric divisors of $\widetilde{\mathbb{P}}_{\!\Delta}$ with $\rho_\Delta=d$. 
Denote by $D_i$ the restriction divisor $r\widetilde{D}_i$ of $\widetilde{D}_i$. 
As to the intersection numbers, 
\begin{equation}\label{SelfIntersection}
D_i^2 = 
\begin{cases}
2l^*(\psi_i)-2 & \textnormal{if $v_i$ is a vertex. }\\
-2 & \textnormal{if $v_i$ is in the interior of an edge. }
\end{cases}
\end{equation}
Here $\psi_i$ is the face in $\Delta^*$ that is dual to $v_i$, and $l^*(\psi_i)$ is the number of lattice points in the interior of the face $\psi_i$. 
Moreover,  
\begin{equation}\label{Intersection}
D_i.D_j = 
\begin{cases}
1 & \textnormal{if $v_i$ and $v_j$ are next to each other on an edge. } \\
l^*(m_{ij}^*)+1 & \textnormal{if $v_i$ and $v_j$ are vertices that are connected } \\
 & \textnormal{by an edge $m_{ij}$ whose dual is $m_{ij}^*$. } \\
0 & \textnormal{otherwise}. 
\end{cases}
\end{equation}
For the above formulas $(2), (3)$, see \cite{FultonToric} (in particular Chapter $5$). 

Recall the following corollaries on lattices by Nikulin~\cite{Nikulin80} and Nishiyama's result~\cite{Nishiyama96} that is needed in the next section. 
\begin{cor}[Corollary 1.6.2~\cite{Nikulin80}]\label{orthogonal}
Let lattices $S$ and $T$ be primitively embedded into the $K3$ lattice. 
Then $S$ and $T$ are orthogonal to each other in the $K3$ lattice if and only if $q_S\simeq -q_T$, where $q_S$ (resp. $q_T$) is the discriminant form of $S$ (resp. $T$). 
\end{cor}

\begin{cor}[Corollary 1.12.3~\cite{Nikulin80}]\label{primitive}
Let $S$ be an even lattice of signature $(t_+, t_-)$ and $\Lambda$ be an even unimodular lattice of signature $(l_+, l_-)$. 
There exists a primitive embedding of $S$ into $\Lambda$ if and only if the following three conditions are simultaneously satisfied. 
\begin{enumerate}
\item[$(1)$] $l_+-l_- \equiv 0 \mod 8$, 
\item[$(2)$] $l_- -t_- \geq 0$ and $l_+ -t_+ \geq 0$, and 
\item[$(3)$] $\rk \Lambda - \rk S > l(A_S)$. 
\end{enumerate}
Here $A_S$ denotes the discriminant group of $S$, which is finitely-generated abelian, and $l(A_S)$ is the minimal number of generators of $A_S$. 
\end{cor}

\begin{remark}
Note that the $K3$ lattice $\Lambda_{K3}$ is of signature $(l_+,l_-)=(3,19)$, so that $l_+-l_- = 19-3=16\equiv 0\mod 8$. 
Thus in order to show that an even hyperbolic lattice $S$ is primitively embedded into $\Lambda_{K3}$, it suffices to verify the conditions $(2)$ and $(3)$ in Corollary \ref{primitive}.  
\end{remark}

\begin{lem}[Lemma 4.3 \cite{Nishiyama96}]\label{NLem4.3}
There exist primitive embeddings of lattices of type $A_1,\, A_2,$ and $A_3$ into $E_8$ with orthogonal complements being respectively $E_7,\, E_6,$ and $D_5$. 
All the notation follows Bourbaki. 
\end{lem}

\section{Main result} \label{MainThm}
We state the main result and prove it in this section. 
\begin{thm}\label{MainThmPf}
With the reflexive polytopes $\Delta$ and $\Delta'$ obtained in \cite{Mase16I}, the families $\mathcal{F}_\Delta$ and $\mathcal{F}_{\Delta'}$ associated to the pairs $(B,\, B')$ of singularities 
\[
(Z_{1,0},\, Z_{1,0}),\, 
(U_{1,0},\, U_{1,0}),\, 
(Q_{17},\, Z_{2,0}),\,  
(W_{1,0},\, W_{1,0})
\]
are lattice dual. 
\end{thm}

Moreover, there exist negative-definite even lattices $K_1$ and $K_2$ such that $\rk{K_1}=15$ and $\discr{K_1}=-18$, and $\rk{K_2}=16$ and $\discr{K_2}=4$ as in Table \ref{LatticeExcPairs}. 
\[
\begin{array}{c @{\hspace{.5cm}} c @{\hspace{.5cm}} c @{\hspace{1cm}} c @{\hspace{1cm}} c @{\hspace{.5cm}} c @{\hspace{.5cm}} c} 
B & \Pic_\Delta & \rho_\Delta & |\discr| & \rho_{\Delta'} & \Pic_{\Delta'} & B' \\ 
\hline \hline
Z_{1,0} & U\oplus D_5\oplus E_7  & 14 & 8 & 6 & U\oplus A_1\oplus A_3 & Z_{1,0}\\
\hline
U_{1,0} & U\oplus K_1& 17  & 18 & 3 & \left( \begin{smallmatrix}0 & 3 \\ 3 & -2\end{smallmatrix}\right)\oplus A_1 & U_{1,0}\\
\hline
Q_{17} & U\oplus E_6\oplus E_7 & 15 & 6 & 5 & U\oplus A_1\oplus A_2 & Z_{2,0}  \\
\hline
W_{1,0} & U\oplus K_2 & 18 & 4 & 2 & \left( \begin{smallmatrix} 0 & 2 \\ 2 & -2 \end{smallmatrix}\right) & W_{1,0}\\
\hline
\end{array}
\]
\begingroup
\captionof{table}{Picard lattices for lattice-dual pairs}\label{LatticeExcPairs}
\endgroup

\noindent
{\sc Proof.}  The statement is verified by an explicit computer-aided calculation case by case. 
In the following, for a lattice $L$ with a basis $\mathcal{B}$, denote by $L^*:={\rm Hom}(L,\mathbb{Z})$ the dual lattice with the dual basis $\mathcal{B}^*$, and $A_L=L^*\slash L$ the discriminant group of $L$, and $l(A_L)$ the minimal number of generators of $A_L$.
In the following, we always set $(l_+, l_-):=(3,19)$ that is the signature of the $K3$ lattice $\Lambda_{K3}$. 
\\
{\bf \mathversion{bold}$Z_{1,0}$ case. } 
Let $\Delta$ and $\Delta'$ be the reflexive polytopes obtained in \cite{Mase16I}. 

Primitive one-simplices of $\Delta'$ are vectors
\[
\begin{array}{llll}
v_1=(-1,0,1) & v_2=(-1,0,0) & v_3=(0,1,-1) & v_4=(2,3,-1) \\
v_5=(2,2,-1) & v_6=(1,-1,-1) & v_7=(0,-1,-1) & v_8=(0,0,-1)\\
v_9=(1,2,-1). 
\end{array}
\]
Denote by $\widetilde{D}_i'$ the the toric divisor corresponding to a one-simplex $v_i$, and $D_i':=r\widetilde{D}_i'$ be the restricted divisor. 
By solving the linear system $(\ref{ToricLinearRelation})$, there is a lattice $L'$ with a basis $\mathcal{B}':=\{ D_i'\, |\, i\not=1,5,6\}$. 
The intersection matrix $M'$ of $L'$ with respect to $\mathcal{B}'$ is computed by formulas (\ref{SelfIntersection}) and (\ref{Intersection}) : 
\[
M'=
\begin{pmatrix}
 -2 & 1 & 0 & 1 & 0 & 0 \\
 1 & -2 & 0 & 0 & 1 & 1 \\
 0 & 0 & -2 & 0 & 0 & 1 \\
 1 & 0 & 0 & -2 & 1 & 0 \\
 0 & 1 & 0 & 1 & -2 & 0 \\
 0 & 1 & 1 & 0 & 0 & -2
\end{pmatrix}
. 
\]
By taking a new basis 
\[
\left\{
\begin{array}{l}
 2 D'_2+5 D'_3+D'_4+2 D'_7+4 D'_8+2 D'_9,\\
 2 D'_2+5 D'_3+D'_4+2 D'_7+3 D'_8+2 D'_9, \\
 3 D'_2+7 D'_3+D'_4+3 D'_7+5 D'_8+2 D'_9, \\
 -D'_2-4 D'_3-D'_4-2 D'_7-3 D'_8-2 D'_9, \\
 -D'_2-2 D'_3-D'_4-D'_8-D'_9, \\
 2 D'_2+6 D'_3+D'_4+2 D'_7+4 D'_8+3 D'_9
\end{array}
\right\},
\]
the lattice $L'$ is isometric to $U\oplus A_1\oplus A_3$. 
It is well-known that $U\oplus A_1\oplus A_3$ is primitively embedded into the $K3$ lattice, and that $A_{L'}\simeq \C{2}\oplus\C{4}$. 
Therefore, $L'=\Pic_{\Delta'}$. 
By a direct computation, $\rk \Pic_{\Delta'}=6$ and $\discr{\Pic_{\Delta'}}=-8$. 
\\

Primitive one-simplices of $\Delta$ are vectors
\[
\begin{array}{llll}
m_1=(1,-1,0) & m_2=(-2,0,-3) & m_3=(-6,2,-7) & m_4=(0,2,-1) \\ 
m_5=(1,1,0) & m_6=(1,0,1) & m_7=(0,0,1) & m_8=(1,0,0) \\
m_9=(-4,1,-5) & m_{10}=(-1,0,-1) & m_{11}=(-3,1,-3) & m_{12}=(0,1,0) \\
m_{13}=(-5,2,-6) & m_{14}=(-4,2,-5) & m_{15}=(-3,2,-4) & m_{16}=(-2,2,-3) \\ 
m_{17}=(-1,2,-2). 
\end{array}
\]
Denote by $\widetilde{D}_i$ the toric divisor corresponding to a one-simplex $m_i$, and $D_i:=r\widetilde{D}_i$ be the restriction divisor. 
By solving the linear system $(\ref{ToricLinearRelation})$, there is a lattice $L$ with a basis $\mathcal{B}=\{ e_k\, | \, k=1,\ldots , 14\} := \{ D_i \, | \, i\not=1,7,10\}$. 
The intersection matrix $M$ of $L$ with respect to $\mathcal{B}$ is computed by formulas (\ref{SelfIntersection}) and (\ref{Intersection}) :
\[
M=
\begin{pmatrix}
 -2 & 0 & 0 & 0 & 0 & 0 & 1 & 0 & 0 & 0 & 0 & 0 & 0 & 0 \\
 0 & -2 & 0 & 0 & 0 & 0 & 1 & 1 & 0 & 1 & 0 & 0 & 0 & 0 \\
 0 & 0 & -2 & 1 & 0 & 0 & 0 & 0 & 1 & 0 & 0 & 0 & 0 & 1 \\
 0 & 0 & 1 & -2 & 1 & 1 & 0 & 0 & 0 & 0 & 0 & 0 & 0 & 0 \\
 0 & 0 & 0 & 1 & -2 & 0 & 0 & 0 & 0 & 0 & 0 & 0 & 0 & 0 \\
 0 & 0 & 0 & 1 & 0 & -2 & 0 & 0 & 0 & 0 & 0 & 0 & 0 & 0 \\
 1 & 1 & 0 & 0 & 0 & 0 & -2 & 0 & 0 & 0 & 0 & 0 & 0 & 0 \\
 0 & 1 & 0 & 0 & 0 & 0 & 0 & -2 & 0 & 0 & 0 & 0 & 0 & 0 \\
 0 & 0 & 1 & 0 & 0 & 0 & 0 & 0 & -2 & 0 & 0 & 0 & 0 & 0 \\
 0 & 1 & 0 & 0 & 0 & 0 & 0 & 0 & 0 & -2 & 1 & 0 & 0 & 0 \\
 0 & 0 & 0 & 0 & 0 & 0 & 0 & 0 & 0 & 1 & -2 & 1 & 0 & 0 \\
 0 & 0 & 0 & 0 & 0 & 0 & 0 & 0 & 0 & 0 & 1 & -2 & 1 & 0 \\
 0 & 0 & 0 & 0 & 0 & 0 & 0 & 0 & 0 & 0 & 0 & 1 & -2 & 1 \\
 0 & 0 & 1 & 0 & 0 & 0 & 0 & 0 & 0 & 0 & 0 & 0 & 1 & -2
\end{pmatrix}. 
\]

In order to prove that the lattice $L$ is primitively embedded into the $K3$ lattice by using Corollary \ref{primitive}, we claim the following statements: \\
{\it Claim.} (1)\, $\rk L =  14, \discr L  =  -8, \sign L = (t_+, t_-) =  (1,13)$. \\
(2)\, $l(A_L)  =  2$. \\
{\it Proof of Claim. }(1) is verified by a direct computation. \\
(2) \, The Gram matrix of $L^*$ with respect to $\mathcal{B}^*:=\{ e_k^*\, |\, k=1,\ldots , 14\}$  is $M^{-1}$. 
We claim that there exist elements of order $4$ and of order $2$ in $A_L$. 
Define a map $\Phi\colon\mathbb{Z}^{15}\to\mathbb{Q}^{14}$ by 
\[
\begin{array}{l}
\Phi[x_1, \, x_2, \, x_3, \, x_4, \, x_5, \, x_6, \, x_7, \, x_8, \, x_9,  \, x_{10}, \, x_{11}, \, x_{12}, \, x_{13}, \, x_{14}, \, s] := \\
\hspace{28mm} s (x_1, \, x_2, \, x_3, \, x_4, \, x_5, \, x_6, \, x_7, \, x_8, \, x_9,  \, x_{10}, \, x_{11}, \, x_{12}, \, x_{13}, \, x_{14}).M^{-1}.
\end{array}
\]
We have 
\[
\Phi[x_1, x_2, x_3, x_4, x_5, x_6, x_7, x_8, x_9, x_{10}, x_{11}, x_{12}, x_{13}, x_{14}; 4]\in\mathbb{Z}^{14}
\]
Thus, for all $\mathbb{X}\in\mathbb{Z}^{14}$, the element $\mathbb{X}.M^{-1}\in\mathbb{Q}^{14}$ satisfies  $4\mathbb{X}.M^{-1}\in\mathbb{Z}^{14} $. 

When $x_5=1$ and $x_6=0$, we have 
\[
\Phi[x_1, x_2, x_3, x_4, 1, 0, x_7, x_8, x_9, x_{10}, x_{11}, x_{12}, x_{13}, x_{14}; 2]\in\mathbb{Q}^{14}\backslash\mathbb{Z}^{14}
\]
Thus by taking general coefficients $(x_i)_{i=1}^{14}$ we get an element of order $4$  
\[
\mathbf{x} = \sum_{i\not=5,6}x_iD_i + D_5
\]
in $A_L$. 
Besides, when $x_5=x_6=0$, we have
\[
\Phi[x_1, x_2, x_3, x_4, 0 , 0, x_7, x_8, x_9, x_{10}, x_{11}, x_{12}, x_{13}, x_{14}; 2] \in\mathbb{Z}^{14}. 
\]
Thus by taking general coefficients $(x_i)_{i=1}^{14}$ we get an element of order $2$  
\[
\mathbf{y} = \sum_{i\not=5,6}x_iD_i
\]
in $A_L$. 
Therefore, 
\[
A_L\simeq \C{4}\oplus \C{2}. 
\]
This completes the proof of the Claim. 
Now we apply Corollary \ref{primitive} to prove that the lattice $L$ is primitively embedded into the $K3$ lattice: 
One obtains
\begin{eqnarray*}
l_+ -t_+  = 3-1=2\geq 0, & & 
l_- -t_-  = 19-12 = 7\geq 0,\\
\rk\Lambda_{K3}-\rk L  =  9 & > & 
l(A_L) = 2.
\end{eqnarray*}
Thus by Corollary \ref{primitive}, the lattice $L$ is primitively embedded into the $K3$ lattice, thus $L=\Pic_\Delta$. 
Besides, since
\[
A_{\Pic_{\Delta}\oplus \, U}\simeq A_{\Pic_{\Delta'}}\simeq\C{2} \oplus \C{4} \textnormal{\, and }\discr (\Pic_{\Delta}\oplus \, U)=-\discr \Pic_{\Delta'}=8, 
\]
by Corollary \ref{orthogonal}, $\Pic_{\Delta}\oplus \, U$ is orthogonal to $\Pic_{\Delta'}$. 

By consulting a list in Lemma \ref{NLem4.3}~\cite{Nishiyama96}, we obtain $\Pic_\Delta\simeq U\oplus E_6 \oplus E_8$. 
\\

\noindent
{\bf \mathversion{bold}$U_{1,0}$ case.} 
Let $\Delta$ and $\Delta'$ be the reflexive polytopes obtained in \cite{Mase16I}. 

Primitive one-simplices of $\Delta'$ are vectors
\[
\begin{array}{llll}
v_1=(-1,0,2) & v_2=(0,1,0) & v_3=(1,2,-1) & v_4=(1,1,-1) \\ 
v_5=(0,-1,0) & v_6=(0,-1,-1). 
\end{array}
\]
Denote by $\widetilde{D}_i'$ the toric divisor corresponding to a one-simplex $v_i$, and $D_i':=r\widetilde{D}_i'$ be the restricted divisor. 
By solving the linear system $(\ref{ToricLinearRelation})$, there is a lattice $L'$ with a basis $\mathcal{B}':=\{ D_2',\, D_3',\, D_4'\}$. 
The intersection matrix $M'$ of $L'$ with respect to $\mathcal{B}'$ is computed by formulas (\ref{SelfIntersection}) and (\ref{Intersection}) : 
\[
M' = 
\begin{pmatrix}
 -2 & 2 & 0 \\
 2 & -2 & 3 \\
 0 & 3 & -2
\end{pmatrix}. 
\]
By taking a new basis 
\[
\left\{
\begin{array}{lll}
D_2'+D_3', & D_4', &  -D_2'
\end{array}
\right\}, 
\]
the lattice $L'$ is isometric to $\left(\mathbb{Z}^3,\, \left( \begin{smallmatrix} 0 & 3 \\ 3 & -2 \end{smallmatrix}\right)\oplus A_1\right)$. \\

\noindent
In order to prove that the lattice $L'$ is primitively embedded into the $K3$ lattice by using Corollary \ref{primitive}, we claim the following statements: \\
{\it Claim 1.} (1)\, $\rk L' = 3, \discr L'  =  18, \sign L'   = (t_+, t_-)  =  (1,2)$. \\
(2)\, $l(A_{L'})  =  2$. \\
{\it Proof of Claim 1. }(1) is verified by a direct computation. \\
(2) \, The Gram matrix of ${L'}^*$ with respect to $\mathcal{B}'^*:=\{ {e'_k}^*\, |\, k=1,2,3\}$  is ${M'}^{-1}$. 
We claim that there exists at least one element of order $9$ in $A_{L'}$. 
Define a map $\Phi\colon\mathbb{Z}^{4}\to\mathbb{Q}^{3}$ by 
\[
\Phi[x_1, \, x_2, \, x_3;\, s] := s (x_1, \, x_2, \, x_3).M'^{-1}.
\]
We have 
\[
\Phi[2x_1, x_2, x_3; 9] \in\mathbb{Z}^3
\]
and 
\[
\Phi[2x_1, x_2, x_3; 3] \in\mathbb{Q}^3\backslash\mathbb{Z}^3. 
\]
Thus by taking general coefficients $(x_1,\, x_2,\, x_3)$ with $x_1\not=0$ or $x_2\not=0$, we get an element of order $9$ 
\[
\mathbf{x}' = 2 x_1D'_1 + x_2D'_2 + x_3D'_3
\]
in $A_{L'}$. 
Besides, we have
\[
\Phi[-1,1, 3x_3; 2] \in\mathbb{Z}^3
\]
Thus we get an element of order $2$  
\[
\mathbf{y}' = -D'_1 + D'_2 + 3x_3D'_3
\]
in $A_{L'}$. 
Therefore, $A_{L'}\simeq \C{9}\oplus\C{2}$. 
In particular, $l(A_{L'})=2$. 
This completes the proof of Claim 1. 
Now we apply Corollary \ref{primitive} to prove that the lattice $L'$ is primitively embedded into the $K3$ lattice: 
One obtains
\begin{eqnarray*}
l_+ -t_+  = 3-1=2\geq 0, & & 
l_- -t_-  = 19-2 = 17\geq 0,\\
\rk\Lambda_{K3}-\rk L'  =  19 & > & 
l(A_{L'}) = 2 .
\end{eqnarray*}
Thus by Corollary \ref{primitive}, the lattice $L'$ is primitively embedded into the $K3$ lattice, thus $L'=\Pic_{\Delta'}$. 
\\

Primitive one-simplices of $\Delta$ are vectors
\[
\begin{array}{llll}
m_1=(1,-1,0) & m_2=(-3,0,-2) & m_3=(-5,1,-3) & m_4=(1,1,0) \\
m_5=(5,-1,2) & m_6=(3,-1,2) & m_7=(0,0,1) & m_8=(-2,1,0) \\ 
m_9=(3,-1,1) & m_{10}=(2,-1,1) & m_{11}=(3,0,1) & m_{12}=(4,-1,2) \\
m_{13}=(-2,0,-1) & m_{14}=(-1,0,0) & m_{15}=(-4,1,-2) & m_{16}=(-3,1,-1) \\
m_{17}=(-1,1,0) & m_{18}=(0,1,0) & m_{19}=(-3,1,-2) & m_{20}=(-1,1,-1). 
\end{array}
\]
Denote by $\widetilde{D}_i$ the toric divisor corresponding to a one-simplex $m_i$, and $D_i:=r\widetilde{D}_i$ be the restricted divisor. 
By solving the linear system $(\ref{ToricLinearRelation})$, there is a lattice $L$ with a basis $\mathcal{B}=\{ e_k\, |\, k=1,\ldots ,17\} := \{ D_i \, | \, i\not=1,7,10\}$. 
The intersection matrix $M$ of $L$ with respect to $\mathcal{B}$ is computed by formulas (\ref{SelfIntersection}) and (\ref{Intersection}): 
\[
M=
\begin{pmatrix}
 -2 & 1 & 0 & 0 & 0 & 0 & 0 & 0 & 0 & 1 & 0 & 0 & 0 & 0 & 0 & 0 & 0 \\
 1 & -2 & 0 & 0 & 0 & 0 & 0 & 0 & 0 & 0 & 0 & 1 & 0 & 0 & 0 & 1 & 0 \\
 0 & 0 & -2 & 0 & 0 & 0 & 0 & 1 & 0 & 0 & 0 & 0 & 0 & 0 & 1 & 0 & 1 \\
 0 & 0 & 0 & -2 & 0 & 0 & 1 & 1 & 1 & 0 & 0 & 0 & 0 & 0 & 0 & 0 & 0 \\
 0 & 0 & 0 & 0 & -2 & 0 & 0 & 0 & 1 & 0 & 0 & 0 & 0 & 0 & 0 & 0 & 0 \\
 0 & 0 & 0 & 0 & 0 & -2 & 0 & 0 & 0 & 0 & 0 & 0 & 1 & 1 & 0 & 0 & 0 \\
 0 & 0 & 0 & 1 & 0 & 0 & -2 & 0 & 0 & 0 & 0 & 0 & 0 & 0 & 0 & 0 & 0 \\
 0 & 0 & 1 & 1 & 0 & 0 & 0 & -2 & 0 & 0 & 0 & 0 & 0 & 0 & 0 & 0 & 0 \\
 0 & 0 & 0 & 1 & 1 & 0 & 0 & 0 & -2 & 0 & 0 & 0 & 0 & 0 & 0 & 0 & 0 \\
 1 & 0 & 0 & 0 & 0 & 0 & 0 & 0 & 0 & -2 & 1 & 0 & 0 & 0 & 0 & 0 & 0 \\
 0 & 0 & 0 & 0 & 0 & 0 & 0 & 0 & 0 & 1 & -2 & 0 & 0 & 0 & 0 & 0 & 0 \\
 0 & 1 & 0 & 0 & 0 & 0 & 0 & 0 & 0 & 0 & 0 & -2 & 1 & 0 & 0 & 0 & 0 \\
 0 & 0 & 0 & 0 & 0 & 1 & 0 & 0 & 0 & 0 & 0 & 1 & -2 & 0 & 0 & 0 & 0 \\
 0 & 0 & 0 & 0 & 0 & 1 & 0 & 0 & 0 & 0 & 0 & 0 & 0 & -2 & 1 & 0 & 0 \\
 0 & 0 & 1 & 0 & 0 & 0 & 0 & 0 & 0 & 0 & 0 & 0 & 0 & 1 & -2 & 0 & 0 \\
 0 & 1 & 0 & 0 & 0 & 0 & 0 & 0 & 0 & 0 & 0 & 0 & 0 & 0 & 0 & -2 & 1 \\
 0 & 0 & 1 & 0 & 0 & 0 & 0 & 0 & 0 & 0 & 0 & 0 & 0 & 0 & 0 & 1 & -2
\end{pmatrix}. 
\]

\noindent
In order to prove that the lattice $L$ is primitively embedded into the $K3$ lattice by using Corollary \ref{primitive}, we claim the following statements: \\
{\it Claim 2.} (1)\, $\rk L = 17, \discr L  = 18, \sign L  = (t_+, t_-)  =  (1,16)$. \\
(2)\, $l(A_{L})  =  2$. \\
{\it Proof of Claim 2. }(1) is verified by a direct computation. \\
(2)\, The Gram matrix of $L^*$ with respect to $\mathcal{B}^*=\{ e_k^*\, |\, k=1,\ldots , 17\}$ is $M^{-1}$. 
We claim that there exists at least one element of order $9$ in $A_L$. 
Define a map $\Phi\colon\mathbb{Z}^{18}\to\mathbb{Q}^{17}$ by 
\[
\begin{array}{l}
\Phi[x_1, \, x_2, \, x_3, \, x_4, \, x_5, \, x_6, \, x_7, \, x_8, \, x_9,  \, x_{10}, \, x_{11}, \, x_{12}, \, x_{13}, \, x_{14}, \, x_{15},\, x_{16},\, x_{17};\, s] :=  \\
\hspace{28mm} s (x_1, \, x_2, \, x_3, \, x_4, \, x_5, \, x_6, \, x_7, \, x_8, \, x_9,  \, x_{10}, \, x_{11}, \, x_{12}, \, x_{13}, \, x_{14}, \, x_{15},\, x_{16},\, x_{17}).M^{-1}.
\end{array}
\]
We have 
\[
\Phi[2x_1, x_2, x_3, x_4, x_5, 2x_6, 2x_7, 2x_8, x_9, x_{10}, 2x_{11}, 2x_{12}, x_{13}, x_{14}, 2x_{15}, x_{16}, x_{17}; 9] \in\mathbb{Z}^{17}, 
\]
and
\[
\Phi[2x_1, x_2, x_3, x_4, x_5, 2x_6, 2x_7, 2x_8, x_9, x_{10}, 2x_{11}, 2x_{12}, x_{13}, x_{14}, 2x_{15}, x_{16}, x_{17}; 3] \in \mathbb{Q}^{17}\backslash \mathbb{Z}^{17}
\]
Let $I=\{ 1,6,7,8,11,12,15\}, \, J=\{ 2,3,4,5,9,10,13,14,16,17\}$. 
Thus by taking general coefficients $(x_i)_{i=1}^{17}$ with $x_{14}\not= 0$, we get an element of order $9$
\[
\mathbf{x} = 2\sum_{i\in I}x_iD_i + \sum_{j\in J}x_iD_i
\]
in $A_L$. 
Besides, we have
\[
\Phi[x_1, 0, x_3, x_4, x_5, 0, x_7, x_8, x_9, 0, 0, 0, 0, 0, 0, 0, 0, 2]  \in\mathbb{Z}^{17}
\]
Thus by taking general coefficients $(x_i)_{i=1,3,4,5,7,8,9}$ we get an element of order $2$
\[
\mathbf{y} = x_1D_1 + x_3D_3 + x_4D_4 + x_5D_5 + x_7D_7 + x_8D_8 + x_9D_9
\]
in $A_L$. 
Therefore,  $A_L=\langle \mathbf{x},\, \mathbf{y}\rangle$, and 
\[
A_L\simeq \C{9}\oplus \C{2}. 
\]
This completes the proof of Claim 2. 
Now we apply Corollary \ref{primitive} to prove that the lattice $L$ is primitively embedded into the $K3$ lattice: 
One obtains
\begin{eqnarray*}
l_+ -t_+  = 3-1=2\geq 0, & & 
l_- -t_-  = 19-16 = 3\geq 0, \\ 
\rk\Lambda_{K3}-\rk L  =  5 &> & 
l(A_L) =2.
\end{eqnarray*}
Thus by Corollary \ref{primitive}, the lattice $L$ is primitively embedded into the $K3$ lattice, thus, $L=\Pic_\Delta$. 
Besides, since
\[
A_{\Pic_{\Delta}\oplus \, U}\simeq A_{\Pic_{\Delta'}}\simeq \C{9}\oplus \C{2}, \textnormal{\, and }\discr(\Pic_{\Delta}\oplus \, U)=-\discr\Pic_{\Delta'}=-18, 
\]
by Corollary \ref{orthogonal}, $\Pic_{\Delta}\oplus \, U$ is orthogonal to $\Pic_{\Delta'}$. 

Since $\rk\Pic_\Delta=17>12$, the lattice $\Pic_\Delta$ contains a sublattice $U$. 
Therefore, there exists a negative-definite even lattice $K_1$ of $\rk{K_1}=17-2=15$ and $\discr{K_1}=-18$. 
\\

\noindent 
{\bf \mathversion{bold}$Q_{17}$ and $Z_{2,0}$ case. }
We first note that all the abelian groups of order $6$ are isomorphic to $\C{6}$. 

Let $\Delta$ and $\Delta'$ be the reflexive polytopes obtained in \cite{Mase16I}. 

Primitive one-simplices of $\Delta'$ are vectors
\[
\begin{array}{llll}
v_1=(-1,-1,2) & v_2=(0,-1,0) & v_3=(1,-1,0) & v_4=(1,-1,1) \\
v_5=(1,2,-3) & v_6=(0,0,-1) & v_7=(1,0,-1) & v_8=(1,1,-2). 
\end{array}
\]
Denote by $\widetilde{D}_i'$ the toric divisor corresponding to a one-simplex $v_i$, and $D_i':=r\widetilde{D}_i'$ be the restricted divisor. 
By solving the linear system $(\ref{ToricLinearRelation})$, there is a lattice $L'$ with a basis $\mathcal{B}':=\{ D_i' \, |\, i\not=1,4,5\}$. 
The intersection matrix $M'$ of $L'$ with respect to $\mathcal{B}'$ is computed by formulas (\ref{SelfIntersection}) and (\ref{Intersection}) :
\[
M'=
\begin{pmatrix}
 -2 & 1 & 2 & 0 & 0 \\
 1 & -2 & 0 & 1 & 0 \\
 2 & 0 & -2 & 0 & 0 \\
 0 & 1 & 0 & -2 & 1 \\
 0 & 0 & 0 & 1 & -2
\end{pmatrix}. 
\]
By taking a new basis 
\[
\left\{
\begin{array}{lllll}
 D'_2+D'_6, & D'_2+D'_3+D'_6, & -D'_6, & -D'_2-D'_6+D'_7, & D'_8
\end{array}
\right\}, 
\]
the lattice $L'$ is isometric to $U\oplus A_1\oplus A_2$. 
It is well-known that $U\oplus A_1\oplus A_2$ is primitively embedded into the $K3$ lattice, and that $A_{L'}\simeq\C{6}$. 
Therefore, $L'=\Pic_{\Delta'}$. 
By a direct computation, $\rk\Pic_{\Delta'}=5$ and $\discr\Pic_{\Delta'}=6$. 
\\

Primitive one-simplices of $\Delta$ are vectors
\[
\begin{array}{llll}
m_1=(0,1,0) & m_2=(-1,0,0) & m_3=(-3,-8,-6) & m_4=(4,-1,1) \\
m_5=(2,1,1) & m_6=(0,1,1) & m_7=(3,0,1) & m_8=(1,1,1) \\ 
m_9=(2,0,1) & m_{10}=(-2,-4,-3) & m_{11}=(-1,-2,-2) & m_{12}=(-2,-5,-4) \\ 
m_{13}=(-2,-7,-5) & m_{14}=(-1,-6,-4) & m_{15}=(0,-5,-3) & m_{16}=(1,-4,-2) \\ 
m_{17}=(2,-3,-1) & m_{18}=(3,-2,0). 
\end{array}
\]
Denote by $\widetilde{D}_i$ the toric divisor corresponding to a one-simplex $m_i$, and $D_i:=r\widetilde{D}_i$ be the restricted divisor. 
By solving the linear system $(\ref{ToricLinearRelation})$, there is a lattice $L$ with a basis $\mathcal{B}=\{ e_k\, |\, e=1,\ldots ,15\} := \{ D_i \, | \, i\not=1,2,6\}$. 
The intersection matrix $M$ of $L$ with respect to $\mathcal{B}$ is computed by formulas (\ref{SelfIntersection}) and (\ref{Intersection}): 
\[
M=
\begin{pmatrix}
 -2 & 0 & 0 & 0 & 0 & 0 & 1 & 0 & 1 & 1 & 0 & 0 & 0 & 0 & 0 \\
 0 & -2 & 0 & 1 & 0 & 1 & 0 & 0 & 0 & 0 & 0 & 0 & 0 & 0 & 1 \\
 0 & 0 & -2 & 1 & 1 & 0 & 0 & 0 & 0 & 0 & 0 & 0 & 0 & 0 & 0 \\
 0 & 1 & 1 & -2 & 0 & 0 & 0 & 0 & 0 & 0 & 0 & 0 & 0 & 0 & 0 \\
 0 & 0 & 1 & 0 & -2 & 0 & 0 & 0 & 0 & 0 & 0 & 0 & 0 & 0 & 0 \\
 0 & 1 & 0 & 0 & 0 & -2 & 0 & 0 & 0 & 0 & 0 & 0 & 0 & 0 & 0 \\
 1 & 0 & 0 & 0 & 0 & 0 & -2 & 0 & 0 & 0 & 0 & 0 & 0 & 0 & 0 \\
 0 & 0 & 0 & 0 & 0 & 0 & 0 & -2 & 1 & 0 & 0 & 0 & 0 & 0 & 0 \\
 1 & 0 & 0 & 0 & 0 & 0 & 0 & 1 & -2 & 0 & 0 & 0 & 0 & 0 & 0 \\
 1 & 0 & 0 & 0 & 0 & 0 & 0 & 0 & 0 & -2 & 1 & 0 & 0 & 0 & 0 \\
 0 & 0 & 0 & 0 & 0 & 0 & 0 & 0 & 0 & 1 & -2 & 1 & 0 & 0 & 0 \\
 0 & 0 & 0 & 0 & 0 & 0 & 0 & 0 & 0 & 0 & 1 & -2 & 1 & 0 & 0 \\
 0 & 0 & 0 & 0 & 0 & 0 & 0 & 0 & 0 & 0 & 0 & 1 & -2 & 1 & 0 \\
 0 & 0 & 0 & 0 & 0 & 0 & 0 & 0 & 0 & 0 & 0 & 0 & 1 & -2 & 1 \\
 0 & 1 & 0 & 0 & 0 & 0 & 0 & 0 & 0 & 0 & 0 & 0 & 0 & 1 & -2
\end{pmatrix}. 
\]
In order to prove that the lattice $L$ is primitively embedded into the $K3$ lattice by using Corollary \ref{primitive}, we claim the following statements: \\
{\it Claim.} (1)\, $\rk L = 15, \discr L  = 6, \sign L  = (t_+, t_-)  =  (1,14)$. \\
(2)\, $l(A_{L})  =  1$. \\
{\it Proof of Claim. }(1) is verified by a direct computation. \\
(2)\, Since $\discr{L}=6,\, A_{L}\simeq \C{6}$. 
This completes the proof of Claim. 
Now we apply Corollary \ref{primitive} to prove that the lattice $L$ is primitively embedded into the $K3$ lattice: 
One obtains
\begin{eqnarray*}
l_+ -t_+  = 3-1=2\geq 0, & & 
l_- -t_-  = 19-14 = 5\geq 0,\\
\rk\Lambda_{K3}-\rk L  =  7 & > & 
l(A_{L})=1.
\end{eqnarray*}
Thus by Corollary \ref{primitive}, the lattice $L$ is primitively embedded into the $K3$ lattice, thus $L=\Pic_\Delta$. 
Besides, since
\[
A_{\Pic_{\Delta}\oplus \, U}\simeq A_{\Pic_{\Delta'}}\simeq \C{6}\textnormal{\, and }\discr (\Pic_{\Delta}\oplus \, U)=-\discr\Pic_{\Delta'}=-6, 
\]
by Corollary \ref{orthogonal}, $\Pic_{\Delta}\oplus \, U$ is orthogonal to $\Pic_{\Delta'}$. 

By consulting a list in Lemma \ref{NLem4.3}~\cite{Nishiyama96}, we obtain $\Pic_\Delta\simeq U\oplus E_6\oplus E_7$. 
\\

\noindent
{\bf \mathversion{bold}$W_{1,0}$ case. }
Let $\Delta$ and $\Delta'$ be the reflexive polytopes obtained in \cite{Mase16I}. 

Primitive one-simplices of $\Delta'$ are vectors
\[
\begin{array}{lllll}
v_1=(-1,0,1) & v_2=(-1,0,0) & v_3=(1,2,-1) & v_4=(2,3,-1) & v_5=(0,-1,0). 
\end{array}
\]
Denote by $\widetilde{D}_i'$ the toric divisor corresponding to a one-simplex $v_i$, and $D\i':=r\widetilde{D}_i'$ be the restriced divisor. 
By solving the linear system $(\ref{ToricLinearRelation})$, there is a lattice $L'$ with a basis $\mathcal{B}':=\{ D_2',\, D_3'\}$. 
The intersection matrix $M'$ of $L'$ with respect to $\mathcal{B}'$ is computed by formulas (\ref{SelfIntersection}) and (\ref{Intersection}) :
\[
M' = 
\left(
\begin{matrix}
 0 & 2 \\
 2 & -2
\end{matrix}\right). 
\]
In order to prove that the lattice $L'$ is primitively embedded into the $K3$ lattice by using Corollary \ref{primitive}, we claim the following statements: \\
{\it Claim 1.} (1)\, $\rk L' = 2, \discr L'  = -4, \sign L'  = (t_+, t_-)  =  (1,1)$. \\
(2)\, $l(A_{L'})  = 2$. \\
{\it Proof of Claim 1. }(1) is verified by a direct computation. \\
(2)\, The Gram matrix of $L'^*$ with respect to $\mathcal{B}'^*:=\{ e'^*_1,\, e'^*_2\}$ is ${M'}^{-1}$. 
By a direct computation, for all $x_1,\, x_2\in\mathbb{Z}$, one gets
\[
2(x_1, x_2)M'^{-1} \in\mathbb{Z}^2, 
\]
which means that elements $2e^*_1,\, 2e^*_2\in L'^*\, (x_2,x_3\in\mathbb{Z}$) are in fact  elements of $L'$.  
Therefore, the discriminant group $A_{L'}$ of $L'$ is isomorphic to $(\C{2})^{\oplus 2}$. 
In particular, $l(A_{L'}) =  2$. 
This completes the proof of Claim 1.
Now we apply Corollary \ref{primitive} to prove that the lattice $L'$ is primitively embedded into the $K3$ lattice: 
One obtains
\begin{eqnarray*}
l_+ -t_+  = 3-1=2\geq 0, & & 
l_- -t_-  = 19-1 = 18\geq 0,\\
\rk\Lambda_{K3}-\rk L'  = 20 & > & 
l(A_{L})=2.
\end{eqnarray*}
Thus by Corollary \ref{primitive}, the lattice $L'$ is primitively embedded into the $K3$ lattice, thus $L'=\Pic_{\Delta'}$. \\

Primitive one-simplices of $\Delta$ are vectors
\[
\begin{array}{llll}
m_1=(1,-1,0) & m_2=(-5,1,-6) & m_3=(1,1,0) & m_4=(1,1,4) \\
m_5=(-1,1,2) & m_6=(1,0,0) & m_7=(1,0,2) & m_8=(0,0,1) \\ 
m_9=(-2,0,-3) & m_{10}=(0,1,3) & m_{11}=(-4,1,-4) & m_{12}=(-3,1,-2) \\
m_{13}=(-2,1,0) & m_{14}=(1,1,1) & m_{15}=(1,1,2) & m_{16}=(1,1,3) \\ 
m_{17}=(-4,1,-5) & m_{18}=(-3,1,-4) & m_{19}=(-2,1,-3) & m_{20}=(-1,1,-2) \\
m_{21}=(0,1,-1). 
\end{array}
\]
Denote by $\widetilde{D}_i$ the toric divisor corresponding to a one-simplex $m_i$, and $D_i := \widetilde{D}_i$ be the restriction divisor. 
By solving the linear system $(\ref{ToricLinearRelation})$, there is a lattice $L$ with a basis $\mathcal{B}=\{ e_k\, |\, k=1,\ldots ,18\}:= \{ D_i \, | \, i\not=1,7,10\}$. 
The intersection matrix $M$ of $L$ with respect to $\mathcal{B}$ is computed by formulas (\ref{SelfIntersection}) and (\ref{Intersection}): 
\[
M=
\begin{pmatrix}
 0 & 0 & 0 & 0 & 1 & 1 & 1 & 0 & 0 & 0 & 0 & 0 & 0 & 0 & 0 & 0 & 0 & 0 \\
 0 & -2 & 0 & 0 & 0 & 0 & 1 & 1 & 0 & 0 & 0 & 0 & 0 & 1 & 0 & 0 & 0 & 0 \\
 0 & 0 & -2 & 0 & 1 & 0 & 0 & 0 & 0 & 0 & 1 & 0 & 0 & 0 & 0 & 0 & 0 & 1 \\
 0 & 0 & 0 & -2 & 0 & 1 & 0 & 0 & 0 & 1 & 0 & 0 & 0 & 0 & 0 & 0 & 0 & 0 \\
 1 & 0 & 1 & 0 & -2 & 0 & 0 & 0 & 0 & 0 & 0 & 0 & 0 & 0 & 0 & 0 & 0 & 0 \\
 1 & 0 & 0 & 1 & 0 & -2 & 0 & 0 & 0 & 0 & 0 & 0 & 0 & 0 & 0 & 0 & 0 & 0 \\
 1 & 1 & 0 & 0 & 0 & 0 & -2 & 0 & 0 & 0 & 0 & 0 & 0 & 0 & 0 & 0 & 0 & 0 \\
 0 & 1 & 0 & 0 & 0 & 0 & 0 & -2 & 1 & 0 & 0 & 0 & 0 & 0 & 0 & 0 & 0 & 0 \\
 0 & 0 & 0 & 0 & 0 & 0 & 0 & 1 & -2 & 1 & 0 & 0 & 0 & 0 & 0 & 0 & 0 & 0 \\
 0 & 0 & 0 & 1 & 0 & 0 & 0 & 0 & 1 & -2 & 0 & 0 & 0 & 0 & 0 & 0 & 0 & 0 \\
 0 & 0 & 1 & 0 & 0 & 0 & 0 & 0 & 0 & 0 & -2 & 1 & 0 & 0 & 0 & 0 & 0 & 0 \\
 0 & 0 & 0 & 0 & 0 & 0 & 0 & 0 & 0 & 0 & 1 & -2 & 1 & 0 & 0 & 0 & 0 & 0 \\
 0 & 0 & 0 & 0 & 0 & 0 & 0 & 0 & 0 & 0 & 0 & 1 & -2 & 0 & 0 & 0 & 0 & 0 \\
 0 & 1 & 0 & 0 & 0 & 0 & 0 & 0 & 0 & 0 & 0 & 0 & 0 & -2 & 1 & 0 & 0 & 0 \\
 0 & 0 & 0 & 0 & 0 & 0 & 0 & 0 & 0 & 0 & 0 & 0 & 0 & 1 & -2 & 1 & 0 & 0 \\
 0 & 0 & 0 & 0 & 0 & 0 & 0 & 0 & 0 & 0 & 0 & 0 & 0 & 0 & 1 & -2 & 1 & 0 \\
 0 & 0 & 0 & 0 & 0 & 0 & 0 & 0 & 0 & 0 & 0 & 0 & 0 & 0 & 0 & 1 & -2 & 1 \\
 0 & 0 & 1 & 0 & 0 & 0 & 0 & 0 & 0 & 0 & 0 & 0 & 0 & 0 & 0 & 0 & 1 & -2
\end{pmatrix}. 
\]
In order to prove that the lattice $L$ is primitively embedded into the $K3$ lattice by using Corollary \ref{primitive}, we claim the following statements: \\
{\it Claim 2.} (1)\, $\rk L = 18, \discr L  = -4, \sign L  = (t_+, t_-)  =  (1,17)$. \\
(2)\, $l(A_{L})  = 2$. \\
{\it Proof of Claim 2. }(1) is verified by a direct computation. \\
(2)\, The Gram matrix of $L^*$ with respect to $\mathcal{B}^*:=\{ e^*_k \, |\, k=1,\ldots , 18\}$ is ${M}^{-1}$. 
An element $x = x_1e_1^* + x_2e_2^* + \cdots + x_{18}e_{18}^*\, (x_1,\, x_2,\ldots , x_{18}\in\mathbb{Z}$) is an $s$-torsion, that is, $sx\in L$ if and only if $s(x_1 \quad x_2\quad \cdots \quad x_{18})M^{-1} \in\mathbb{Z}^{18}. $
By a direct computation, for all $x_1,\, x_2,\, \cdots ,\, x_{18}\in\mathbb{Z}$, one gets
\[
2(x_1 \quad x_2\quad \cdots \quad x_{18})M^{-1}\in\mathbb{Z}^{18}. 
\]
Thus, there exist two distinct generating vectors in $A_L$ of order $2$. 
Indeed, one may take 
\[
x = e_1^* + e_3^* + \cdots + e_{17}^*, \quad 
x' = e_2^* + e_4^* + \cdots + e_{18}^* ,
\]
then, $\langle x\rangle = \langle x'\rangle = \mathbb{Z}\slash 2\mathbb{Z},$ and $\langle x\rangle \cap \langle x'\rangle = \{ 0\}$. 
Therefore, 
\[
A_L\simeq (\C{2})^{\oplus 2}. 
\]
This completes the proof of Claim 2. 
Now we apply Corollary \ref{primitive} to prove that the lattice $L$ is primitively embedded into the $K3$ lattice: 
One obtains
\begin{eqnarray*}
l_+ -t_+  = 3-1=2\geq 0, & & 
l_- -t_-  = 19-17 = 2\geq 0 ,\\
\rk\Lambda_{K3}-\rk L  =  4 & > &
l(A_L) = 2.
\end{eqnarray*}
Thus by Corollary \ref{primitive}, the lattice $L$ is primitively embedded into the $K3$ lattice, thus $L=\Pic_\Delta$. 
Besides, since 
\[
A_{\Pic_{\Delta}\oplus \, U}\simeq A_{\Pic_{\Delta'}}\simeq (\C{2})^{\oplus 2} \textnormal{\, and } \discr (\Pic_{\Delta}\oplus \, U)=-\discr\Pic_{\Delta'}=4, 
\]
by Corollary \ref{orthogonal}, $\Pic_{\Delta}\oplus \, U$ is orthogonal to $\Pic_{\Delta'}$. 

Since $\rk\Pic_\Delta=18>12$, the lattice $\Pic_\Delta$ contains a sublattice $U$. 
Therefore, there exists a negative-definite even lattice $K_2$ of $\rk{K_2}=18-2=16$ and $\discr{K_2}=4$. 
Thus the result follows. 
\QED

\begin{remark}
Denote by $\Delta'_B$ for the polytope $\Delta'$ obtained in Theorem \ref{MainThmPf} for the singularity of type $B$ . \\
(1)\, The families of $K3$ surfaces associated to singularities except $\mathcal{F}_{\Delta'_{U_{1,0}}}$ and $\mathcal{F}_{\Delta'_{W_{1,0}}}$ contain general members having the structure of a Jacobian elliptic fibration. \\
(2)\, The family $\mathcal{F}_{\Delta'_{U_{1,0}}}$ of $K3$ surfaces contains general members with the structure of elliptic fibration with a $3$-section. 
Indeed, $\Pic_{\Delta'_{U_{1,0}}}$ contains a sublattice generated by $\langle f, s\rangle $   with respect to which the intersection matrix is $\left( \begin{smallmatrix}0 & 3 \\ 3 & -2\end{smallmatrix}\right)$ that is clearly not isometric to $U$. 
The classes $f$ and $s$ are respectively, of an elliptic curve $E$ and of a nonsingular rational curve $C$ such that $E$ and $C$ intersect at three points, which means that $C$ is a $3$-section. 
\\
(3)\, The family $\mathcal{F}_{\Delta'_{W_{1,0}}}$ of $K3$ surfaces contains general members having the structure of elliptic fibration with a $2$-section. 
Indeed, $\Pic_{\Delta'_{W_{1,0}}}$ contains a sublattice generated by $\langle f, s\rangle $  with respect to which the intersection matrix is $\left( \begin{smallmatrix}0 & 2 \\ 2 & -2\end{smallmatrix}\right)$ that is clearly not isometric to $U$. 
The classes $f$ and $s$ are respectively, of an elliptic curve $E$ and of a nonsingular rational curve $C$ such that $E$ and $C$ intersect at two points, which means that $C$ is a $2$-section. 
\end{remark}

\hfill \textsc{Makiko Mase} \\ 
\hfill e-mail: mtmase@arion.ocn.ne.jp \\
\hfill {\small{Department of Mathematics and Information Sciences, Tokyo Metropolitan University}} \\
\hfill {\small{1-1 Minami-Osawa, Hachioji-shi Tokyo, 192-0397, Japan }}\\
\end{document}